\documentclass[12pt,reqno]{article}

\usepackage[usenames]{color}
\usepackage{graphicx}
\usepackage{amscd}

\usepackage{color}
\usepackage{fullpage}
\usepackage{float}

\usepackage{psfig}
\usepackage{graphics,amsmath,amssymb}
\usepackage{amsfonts}
\usepackage{latexsym}
\usepackage{epsf}
\usepackage{amsthm}
\usepackage{pstricks,pst-node}

\usepackage{xspace}
\usepackage[colorlinks=true,
linkcolor=webgreen,
filecolor=webbrown,
citecolor=webgreen]{hyperref}

\definecolor{webgreen}{rgb}{0,.5,0}
\definecolor{webbrown}{rgb}{.6,0,0}

\def\v{\vert}
\def\p{\ensuremath{\mathcal P}\xspace}
\def\bl{baseline\xspace}
\def\kd{$k$-divisible\xspace}

\def\blueclr{\textcolor{blue} }

\newtheorem*{theorem}{Theorem}
\newtheorem*{main}{Main Theorem}
\newtheorem*{lemma}{Lemma}
\newtheorem*{cor}{Corollary}

\catcode`\@=11

\thicklines
\newskip\Einheit \Einheit=.6cm
\newcount\xcoord \newcount\ycoord
\newdimen\xdim \newdimen\ydim \newdimen\PfadD@cke \newdimen\Pfadd@cke
\def\PfadDicke#1{\PfadD@cke#1 \divide\PfadD@cke by2 
\Pfadd@cke\PfadD@cke \multiply\PfadD@cke by2}
\long\def\LOOP#1\REPEAT{\def\BODY{#1}\ITERATE}
\def\ITERATE{\BODY \let\next\ITERATE \else\let\next\relax\fi \next}
\let\REPEAT=\fi
\def\Punkt{\hbox{\raise-2pt\hbox to0pt{\hss\scriptsize$\bullet$\hss}}}

\def\DuennPunkt(#1,#2){\unskip
  \raise#2 \Einheit\hbox to0pt{\hskip#1 \Einheit
          \raise-1.5pt\hbox to0pt{\hss\tiny$\bullet$\hss}\hss}}
		  
\def\NormalPunkt(#1,#2){\unskip
  \raise#2 \Einheit\hbox to0pt{\hskip#1 \Einheit
          \raise-3pt\hbox to0pt{\hss\large$\bullet$\hss}\hss}}
\def\DickPunkt(#1,#2){\unskip
  \raise#2 \Einheit\hbox to0pt{\hskip#1 \Einheit
          \raise-4pt\hbox to0pt{\hss\Large$\bullet$\hss}\hss}}
\def\Kreis(#1,#2){\unskip
  \raise#2 \Einheit\hbox to0pt{\hskip#1 \Einheit
          \raise-4pt\hbox to0pt{\hss\Large$\circ$\hss}\hss}}
\def\Diagonale(#1,#2)#3{\unskip\leavevmode
  \xcoord#1\relax \ycoord#2\relax
      \raise\ycoord \Einheit\hbox to0pt{\hskip\xcoord \Einheit
         \unitlength\Einheit
         \line(1,1){#3}\hss}}
\def\AntiDiagonale(#1,#2)#3{\unskip\leavevmode
  \xcoord#1\relax \ycoord#2\relax \advance\xcoord by -0.05\relax
      \raise\ycoord \Einheit\hbox to0pt{\hskip\xcoord \Einheit
         \unitlength\Einheit
         \line(1,-1){#3}\hss}}
\def\Pfad(#1,#2),#3\endPfad{\unskip\leavevmode
  \xcoord#1 \ycoord#2 \thicklines\ZeichnePfad#3\endPfad\thinlines}
\def\ZeichnePfad#1{\ifx#1\endPfad\let\next\relax
  \else\let\next\ZeichnePfad
    \ifnum#1=1
      \raise\ycoord \Einheit\hbox to0pt{\hskip\xcoord \Einheit
         \vrule height\Pfadd@cke width1 \Einheit depth\Pfadd@cke\hss}%
      \advance\xcoord by 1
     \else\ifnum#1=2
      \raise\ycoord \Einheit\hbox to0pt{\hskip\xcoord \Einheit
         \unitlength\Einheit
         \line(0,1){1}\hss}
      \advance\xcoord by 0
      \advance\ycoord by 1
 \else\ifnum#1=3
      \raise\ycoord \Einheit\hbox to0pt{\hskip\xcoord \Einheit
         \unitlength\Einheit
         \line(1,1){1}\hss}
      \advance\xcoord by 1
      \advance\ycoord by 1
    \else\ifnum#1=4
      \raise\ycoord \Einheit\hbox to0pt{\hskip\xcoord \Einheit
         \unitlength\Einheit
         \line(1,-1){1}\hss}
      \advance\xcoord by 1
      \advance\ycoord by -1
   \else\ifnum#1=5
      \raise\ycoord \Einheit\hbox to0pt{\hskip\xcoord \Einheit
         \unitlength\Einheit
         \line(2,1){2}\hss}
      \advance\xcoord by 2
      \advance\ycoord by 1
	  \else\ifnum#1=6
      \raise\ycoord \Einheit\hbox to0pt{\hskip\xcoord \Einheit
         \unitlength\Einheit
         \line(2,-1){2}\hss}
      \advance\xcoord by 2
      \advance\ycoord by -1
	  \else\ifnum#1=7
      \raise\ycoord \Einheit\hbox to0pt{\hskip\xcoord \Einheit
         \unitlength\Einheit
         \line(3,1){3}\hss}
      \advance\xcoord by 3
      \advance\ycoord by 1
	  \else\ifnum#1=8
      \raise\ycoord \Einheit\hbox to0pt{\hskip\xcoord \Einheit
         \unitlength\Einheit
         \line(3,-1){3}\hss}
      \advance\xcoord by 3
      \advance\ycoord by -1
    \fi\fi\fi\fi\fi\fi\fi\fi
  \fi\next}
\def\hSSchritt{\leavevmode\raise-.4pt\hbox 
to0pt{\hss.\hss}\hskip.2\Einheit
  \raise-.4pt\hbox to0pt{\hss.\hss}\hskip.2\Einheit
  \raise-.4pt\hbox to0pt{\hss.\hss}\hskip.2\Einheit
  \raise-.4pt\hbox to0pt{\hss.\hss}\hskip.2\Einheit
  \raise-.4pt\hbox to0pt{\hss.\hss}\hskip.2\Einheit}
\def\vSSchritt{\vbox{\baselineskip.2\Einheit\lineskiplimit0pt
\hbox{.}\hbox{.}\hbox{.}\hbox{.}\hbox{.}}}
\def\DSSchritt{\leavevmode\raise-.4pt\hbox to0pt{%
  \hbox to0pt{\hss.\hss}\hskip.2\Einheit
  \raise.2\Einheit\hbox to0pt{\hss.\hss}\hskip.2\Einheit
  \raise.4\Einheit\hbox to0pt{\hss.\hss}\hskip.2\Einheit
  \raise.6\Einheit\hbox to0pt{\hss.\hss}\hskip.2\Einheit
  \raise.8\Einheit\hbox to0pt{\hss.\hss}\hss}}
\def\dSSchritt{\leavevmode\raise-.4pt\hbox to0pt{%
  \hbox to0pt{\hss.\hss}\hskip.2\Einheit
  \raise-.2\Einheit\hbox to0pt{\hss.\hss}\hskip.2\Einheit
  \raise-.4\Einheit\hbox to0pt{\hss.\hss}\hskip.2\Einheit
  \raise-.6\Einheit\hbox to0pt{\hss.\hss}\hskip.2\Einheit
  \raise-.8\Einheit\hbox to0pt{\hss.\hss}\hss}}
\def\SPfad(#1,#2),#3\endSPfad{\unskip\leavevmode
  \xcoord#1 \ycoord#2 \ZeichneSPfad#3\endSPfad}
\def\ZeichneSPfad#1{\ifx#1\endSPfad\let\next\relax
  \else\let\next\ZeichneSPfad
    \ifnum#1=1
      \raise\ycoord \Einheit\hbox to0pt{\hskip\xcoord \Einheit
         \hSSchritt\hss}%
      \advance\xcoord by 1
    \else\ifnum#1=2
      \raise\ycoord \Einheit\hbox to0pt{\hskip\xcoord \Einheit
        \hbox{\hskip-2pt \vSSchritt}\hss}%
      \advance\ycoord by 1
    \else\ifnum#1=3
      \raise\ycoord \Einheit\hbox to0pt{\hskip\xcoord \Einheit
         \DSSchritt\hss}
      \advance\xcoord by 1
      \advance\ycoord by 1
    \else\ifnum#1=4
      \raise\ycoord \Einheit\hbox to0pt{\hskip\xcoord \Einheit
         \dSSchritt\hss}
      \advance\xcoord by 1
      \advance\ycoord by -1
    \fi\fi\fi\fi
  \fi\next}
\def\Koordinatenachsen(#1,#2){\unskip
 \hbox to0pt{\hskip-.5pt\vrule height#2 \Einheit width.5pt depth1 
\Einheit}%
 \hbox to0pt{\hskip-1 \Einheit \xcoord#1 \advance\xcoord by1
    \vrule height0.25pt width\xcoord \Einheit depth0.25pt\hss}}
\def\Koordinatenachsen(#1,#2)(#3,#4){\unskip
 \hbox to0pt{\hskip-.5pt \ycoord-#4 \advance\ycoord by1
    \vrule height#2 \Einheit width.5pt depth\ycoord \Einheit}%
 \hbox to0pt{\hskip-1 \Einheit \hskip#3\Einheit 
    \xcoord#1 \advance\xcoord by1 \advance\xcoord by-#3 
    \vrule height0.25pt width\xcoord \Einheit depth0.25pt\hss}}
\def\Gitter(#1,#2){\unskip \xcoord0 \ycoord0 \leavevmode
  \LOOP\ifnum\ycoord<#2
    \loop\ifnum\xcoord<#1
      \raise\ycoord \Einheit\hbox to0pt{\hskip\xcoord 
\Einheit\Punkt\hss}%
      \advance\xcoord by1
    \repeat
    \xcoord0
    \advance\ycoord by1
  \REPEAT}
\def\Gitter(#1,#2)(#3,#4){\unskip \xcoord#3 \ycoord#4 \leavevmode
  \LOOP\ifnum\ycoord<#2
    \loop\ifnum\xcoord<#1
      \raise\ycoord \Einheit\hbox to0pt{\hskip\xcoord 
\Einheit\Punkt\hss}%
      \advance\xcoord by1
    \repeat
    \xcoord#3
    \advance\ycoord by1
  \REPEAT}
\def\Label#1#2(#3,#4){\unskip \xdim#3 \Einheit \ydim#4 \Einheit
  \def\lo{\advance\xdim by-.5 \Einheit \advance\ydim by.5 \Einheit}%
  \def\llo{\advance\xdim by-.25cm \advance\ydim by.5 \Einheit}%
  \def\loo{\advance\xdim by-.5 \Einheit \advance\ydim by.25cm}%
  \def\o{\advance\ydim by.25cm}%
  \def\ro{\advance\xdim by.5 \Einheit \advance\ydim by.5 \Einheit}%
  \def\rro{\advance\xdim by.25cm \advance\ydim by.5 \Einheit}%
  \def\roo{\advance\xdim by.5 \Einheit \advance\ydim by.25cm}%
  \def\l{\advance\xdim by-.30cm}%
  \def\r{\advance\xdim by.30cm}%
  \def\lu{\advance\xdim by-.5 \Einheit \advance\ydim by-.6 \Einheit}%
  \def\llu{\advance\xdim by-.25cm \advance\ydim by-.6 \Einheit}%
  \def\luu{\advance\xdim by-.5 \Einheit \advance\ydim by-.30cm}%
  \def\u{\advance\ydim by-.30cm}%
  \def\ru{\advance\xdim by.5 \Einheit \advance\ydim by-.6 \Einheit}%
  \def\rru{\advance\xdim by.25cm \advance\ydim by-.6 \Einheit}%
  \def\ruu{\advance\xdim by.5 \Einheit \advance\ydim by-.30cm}%
  #1\raise\ydim\hbox to0pt{\hskip\xdim
     \vbox to0pt{\vss\hbox to0pt{\hss$#2$\hss}\vss}\hss}%
}
\catcode`\@=12

\begin{document}

\begin{center}
\vskip 1cm{\LARGE\bf 
A Combinatorial Interpretation of 
 $\frac{j}{n}\binom{kn}{n+j} $ 

}
\vskip 1cm
\large
David Callan\\
Department of Statistics \\
University of Wisconsin-Madison  \\
1300 University Ave  \\
Madison, WI \ 53706-1532  \\

\href{mailto:callan@stat.wisc.edu}{\tt callan@stat.wisc.edu} \\
\end{center}

\vspace*{5mm}

\begin{abstract}
The identity  $\frac{j}{n}\binom{kn}{n+j} 
=(k-1)\binom{kn-1}{n+j-1}-\binom{kn-1}{n+j}$ shows that  
$\frac{j}{n}\binom{kn}{n+j} $ is always an integer. 
Here we give a combinatorial interpretation of 
this integer in terms of lattice paths, using a 
uniformly distributed statistic. In 
particular, the case $j=1,k=2$ gives yet another manifestation of the 
Catalan numbers.
\end{abstract}

\vspace{6mm}

\section{Introduction} For each pair of integers $j\ge 1$ and $k\ge 2$, the sequence 
$\big(\frac{j}{n}\binom{kn}{n+j}\big)_{n\ge \frac{j}{k-1}}$ consists of integers 
since $\frac{j}{n}\binom{kn}{n+j} 
=(k-1)\binom{kn-1}{n+j-1}-\binom{kn-1}{n+j}$. For $j=1,k=2$ this sequence is 
the Catalan numbers,
\htmladdnormallink{A000108}{http://www.research.att.com:80/cgi-bin/access.cgi/as/njas/sequences/eisA.cgi?Anum=A000108}
in the
\htmladdnormallink{On-Line Encyclopedia of Integer Sequences}{http://www.research.att.com/~njas/sequences/Seis.html}
; for $j=1,k=3$ it is
\htmladdnormallink{A007226}{http://www.research.att.com:80/cgi-bin/access.cgi/as/njas/sequences/eisA.cgi?Anum=A007226}
and for $j=1,k=4$ it is
\htmladdnormallink{A007228}{http://www.research.att.com:80/cgi-bin/access.cgi/as/njas/sequences/eisA.cgi?Anum=A007228}
. In this note, we give a combinatorial interpretation for all 
$j,k$ in terms of lattice paths. We first treat 
the case $j=1$, which is simpler (\S 2), specialize to $k=1$ (\S 3), then generalize 
to larger $j$ (\S 4), and end with some remarks (\S 5). 

\section{Case \emph{j} = 1}
Let $\p_{n,k}$ denote the set of lattice paths of $n+1$ upsteps 
$U=(1,1)$ and $(k-1)n-1$ downsteps $D=(1,-1)$. Clearly, $\v 
\p_{n,k}\v=\binom{kn}{n+1}\,:$ choose locations for the upsteps among 
the total of $kn$ steps. A path in $\p_{n,k}$ has $kn+1$ vertices or 
\emph{points}: its initial and terminal points and $kn-1$ interior 
points. Define the \emph{\bl} of $P\in 
\p_{n,k}$ to be the line joining its initial and terminal points. For $P\in 
\p_{n,k}$ label its points $0,1,2,\ldots,kn$ left to right and 
define the \emph{\kd}points of $P$ to be those whose label is 
divisible by $k$. An example with $n=5$ and $k=3$ is illustrated 
(\kd points indicated by a heavy dot).

\vspace*{-20mm}
\begin{pspicture}(-2,-1)(16,6)
\psset{unit=.7cm}   
\psdots(0,3) 
(1,4)(2,3)(3,2)(4,3)(5,2)(6,3)(7,4)(8,5)(9,4)(10,3)(11,2)(12,1)(13,0)(14,1)(15,0)
\psline[linecolor=blue](0,3)(1,4)(2,3)(3,2)(4,3)(5,2)(6,3)(7,4)(8,5)(9,4)(10,3)(11,2)(12,1)(13,0)(14,1)(15,0)
\psline[linestyle=dotted](0,3)(15,0)

\psset{dotscale=2}
\psdots(0,3)(3,2)(6,3)(9,4)(12,1)(15,0)
\rput(-0.2,2.6){\textrm{{\footnotesize $0$}}}
\rput(1,4.4){\textrm{{\scriptsize $1$}}}
\rput(2.1,3.4){\textrm{{\scriptsize $2$}}}
\rput(3,1.6){\textrm{{\scriptsize $3$}}}
\rput(4,3.4){\textrm{{\scriptsize $4$}}}
\rput(5,1.6){\textrm{{\scriptsize $5$}}}
\rput(5.9,3.4){\textrm{{\scriptsize $6$}}}
\rput(6.9,4.4){\textrm{{\scriptsize $7$}}}
\rput(8,5.4){\textrm{{\scriptsize $8$}}}
\rput(9.1,4.4){\textrm{{\scriptsize $9$}}}
\rput(10.1,3.4){\textrm{{\scriptsize $10$}}}
\rput(11.1,2.4){\textrm{{\scriptsize $11$}}}
\rput(12.1,1.4){\textrm{{\scriptsize $12$}}}
\rput(13,-0.4){\textrm{{\scriptsize $13$}}}
\rput(14,1.4){\textrm{{\scriptsize $14$}}}
\rput(15,-0.4){\textrm{{\scriptsize $15$}}}

\rput(8,0.2){$\textrm{{\footnotesize \bl}}$}
\rput(8,0.9){$\textrm{{\footnotesize $\nearrow$}}$}
\rput(8,-1.3){$\textrm{{\small A path in $\p_{5,3}$}}$}

\end{pspicture}

\vspace*{5mm}

Consider the statistic $X$ on $\p_{n,k}$ defined by $X=\#\:$interior \kd 
 points lying strictly above the \bl. In the illustration 
$X=3$ (points 6,\ 9 and 12).
\begin{theorem}
    The statistic $X$ on $\p_{n,k}$ is uniformly distributed over 
    $0,1,2,\ldots,n-1$.
\end{theorem}
The following count is an immediate consequence of the theorem by  
considering the paths with $X=n-1$.
\begin{cor}
   $\frac{1}{n}\binom{kn}{n+1}$ is the number of paths in $\p_{n,k}$ 
   all of whose interior \kd points lie strictly above its \bl. 
\end{cor}
\noindent \textbf{Proof of Theorem}\quad Consider the operation 
``rotate left $k$ units'' on $\p_{n,k}$ defined by transferring the 
initial $k$ steps of a path in $\p_{n,k}$ to the end. This rotation 
operation partitions $\p_{n,k}$ into rotation classes. We claim (i) 
each such rotation class has size $n$, and (ii) $X$ assumes the values
$0,1,2,\ldots,n-1$ in turn on the paths of a rotation class. The first claim follows 
from
\begin{lemma}
    Given $P \in \p_{n,k}$, the only \kd points lying on its \bl are 
    its initial and terminal points.
\end{lemma}
\noindent \textbf{Proof of Lemma}\quad Suppose $ik,\ 0\le i \le n$, is a \kd 
point on the \bl. Since the slope of the \bl is 
$-\frac{(k-2)n-2}{kn}$, this says that the point with coordinates 
$(ik,-i \frac{(k-2)n-2}{n})$ lies on $P$ (taking the initial point of $P$ as 
origin). For each point $(x,y)$ on $P$, $x$ and $y$ must 
have the same even/odd parity. Hence $ik \equiv i 
\frac{(k-2)n-2}{n}$\:mod\:2. Simplifying, we find $2i \equiv \frac{2i}{n} 
\:\textrm{mod}\:2 \Rightarrow  i \equiv \frac{i}{n} 
\:\textrm{mod}\:1 \Rightarrow n|\,i \Rightarrow i=0$ or $n$, the last 
implication because $0\le i \le n$. \qed

To prove the second claim, we exhibit a bijection from the paths in 
$\p_{n,k}$ with $X=n-1$ to those with $X=i$ for each $i\in [0,n-1]$.
Given $P\in \p_{n,k}$ with $X=n-1$, draw its \bl $L$. The entire 
rotation class of $P$ can be viewed in a single diagram: 
draw a second contiguous copy of $P$ as illustrated, then join the two 
occurrences of each interior \kd point. This results in $n$ parallel 
line segments (no two collinear, by the Lemma), each the base line of 
a path in the rotation class of $P$. 
Label the lines (at their endpoints) 0 through $n-1$ from top to bottom. 

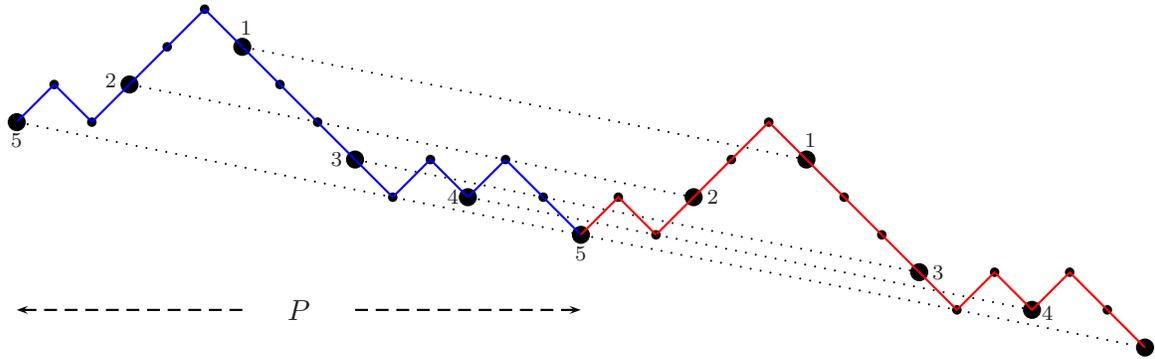
\vspace*{-45mm}
\begin{pspicture}(0,-1.5)(30,10)
\psset{unit=.5cm}   
\psdots(0,6) 
(1,7)(2,6)(3,7)(4,8)(5,9)(6,8)(7,7)(8,6)(9,5)(10,4)(11,5)(12,4)(13,5)(14,4)(15,3)
(16,4)(17,3)(18,4)(19,5)(20,6)(21,5)(22,4)(23,3)(24,2)(25,1)(26,2)(27,1)(28,2)(29,1)(30,0)
\psset{dotscale=1.9}
\psdots(0,6)(3,7)(6,8)(9,5)(12,4)(15,3)
\psdots(18,4)(21,5)(24,2)(27,1)(30,0)
 
\psline[linecolor=blue](0,6) 
(1,7)(2,6)(3,7)(4,8)(5,9)(6,8)(7,7)(8,6)(9,5)(10,4)(11,5)(12,4)(13,5)(14,4)(15,3)

\psline[linecolor=red](15,3)
(16,4)(17,3)(18,4)(19,5)(20,6)(21,5)(22,4)(23,3)(24,2)(25,1)(26,2)(27,1)(28,2)(29,1)(30,0)

\psline[linestyle=dotted](0,6)(30,0)
\psline[linestyle=dotted](3,7)(18,4)
\psline[linestyle=dotted](6,8)(21,5)
\psline[linestyle=dotted](9,5)(24,2)
\psline[linestyle=dotted](12,4)(27,1)

\psline[linestyle=dashed]{<-}(0,1)(6,1)
\psline[linestyle=dashed]{->}(9,1)(15,1)
\rput(7.5,1){\textrm{{\small $P$}}}

\rput(0,5.5){\textrm{{\scriptsize $5$}}}
\rput(15,2.5){\textrm{{\scriptsize $5$}}}
\rput(2.5,7.1){\textrm{{\scriptsize $2$}}}
\rput(18.5,4){\textrm{{\scriptsize $2$}}}
\rput(6.1,8.5){\textrm{{\scriptsize $1$}}}
\rput(21.1,5.5){\textrm{{\scriptsize $1$}}}
\rput(8.5,5){\textrm{{\scriptsize $3$}}}
\rput(24.5,2){\textrm{{\scriptsize $3$}}}
\rput(11.6,4){\textrm{{\scriptsize $4$}}}
\rput(27.4,0.9){\textrm{{\scriptsize $4$}}}

\rput(14,-2){\textrm{{\small The rotation class of 
$P=UDU^{3}D^{5}UDUD^{2} \in \p_{5,3}$}}}

\end{pspicture}


Now the path $Q$ with \bl $i$ has the form $BA$ when $P$ is decomposed 
as $AB$ with $A$ an initial segment of $P$.  Hence $Q$ is in $\p_{n,k}$ and has $X=i$ since the 
interior \kd points of $Q$ lying (strictly) 
above its \bl are precisely those labeled $1,2,\ldots,i-1$. The path $B$ can be 
retrieved in $Q$ as the initial subpath of $Q$ terminating at its 
``lowest'' \kd point where ``lowest'' is measured relative to the 
parallel lines, and so the mapping is invertible.

The diagram used in this proof is reminiscent of the one used in \emph{Concrete 
Mathematics} \cite[p.\,360]{gkp} to  prove Raney's Lemma, also known as the 
Cycle Lemma \cite{zaks,sw}.

\section{Special Case} The case $j=1,k=2$ gives a new interpretation of 
the Catalan numbers: $C_{n}$ is the number of lattice paths of $n+1$ upsteps 
and $n-1$ downsteps such that the interior even-numbered vertices all lie 
strictly 
above the line joining the initial and terminal points. The $C_{3}=5$ 
paths with $n=3$ are shown.
\Einheit=0.4cm
\[
\blueclr{\Pfad(-17,0),333344\endPfad
\Pfad(-10,0),333434\endPfad
\Pfad(-3,0),333443\endPfad
\Pfad(4,0),334334\endPfad
\Pfad(11,0),334343\endPfad}
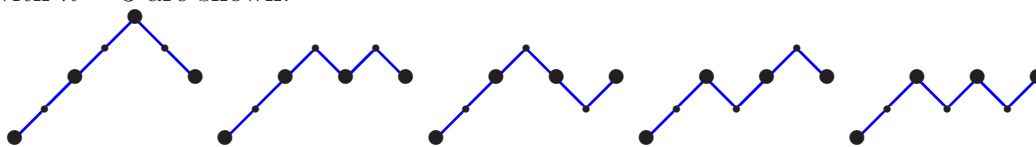
\Label\u{ \textrm{{\scriptsize The 5 paths in $\p_{3,2}$}}}(0,-1)
\NormalPunkt(-17,0)
\DuennPunkt(-16,1)
\NormalPunkt(-15,2)
\DuennPunkt(-14,3)
\NormalPunkt(-13,4)
\DuennPunkt(-12,3)
\NormalPunkt(-11,2)
\NormalPunkt(-10,0)
\DuennPunkt(-9,1)
\NormalPunkt(-8,2)
\DuennPunkt(-7,3)
\NormalPunkt(-6,2)
\DuennPunkt(-5,3)
\NormalPunkt(-4,2)
\NormalPunkt(-3,0)
\DuennPunkt(-2,1)
\NormalPunkt(-1,2)
\DuennPunkt(0,3)
\NormalPunkt(1,2)
\DuennPunkt(2,1)
\NormalPunkt(3,2)
\NormalPunkt(4,0)
\DuennPunkt(5,1)
\NormalPunkt(6,2)
\DuennPunkt(7,1)
\NormalPunkt(8,2)
\DuennPunkt(9,3)
\NormalPunkt(10,2)
\NormalPunkt(11,0)
\DuennPunkt(12,1)
\NormalPunkt(13,2)
\DuennPunkt(14,1)
\NormalPunkt(15,2)
\DuennPunkt(16,1)
\NormalPunkt(17,2)
\]

\vspace*{5mm}

\section{General Case}

The general case $j\ge 1$ is similar but a little more complicated. 
Let $\p_{n,k,j}$ denote the set of paths of $kn$ upsteps/downsteps of 
which $n+j$ are upsteps. Thus $\v \p_{n,k,j} \v =\binom{kn}{n+j}$. 
The ``$j$'' factor in the numerator of 
$\frac{j}{n}\binom{kn}{n+j}$ requires that we consider the Cartesian 
product $\p_{n,k,j}^{*}:= \p_{n,k,j} \times [j]$ whose size is 
$j\binom{kn}{n+j}$. Given $(P,i)\in \p_{n,k,j}^{*}$, introduce an 
$x$-$y$ coordinate system with origin at the initial point of $P$, 
identify the parameter $i$ with the line segment joining $(0,2(i-1))$ 
and $(kn,-(k-2)n+2i)$, and call this the baseline for $(P,i)$; it 
coincides with the previous notion of baseline when $j=1$, forcing 
$i=1$. 
 It is easy to see that, once again, 
the \bl never contains an interior \kd point of $P$. Define $X$ 
on $(P,i)\in \p_{n,k,j}^{*}$ by $X=\#\ $interior \kd points of $P$ 
lying strictly above the \bl.

We first show that $X$ is uniformly distributed over $0,1,2,\ldots,n-1$. 
It is no longer true that orbits in $\p_{n,k,j}$ under the ``rotate 
left by $k$'' operator $R$ all have size $n$ but no matter: in 
general,
$P\in \p_{n,k,j}$ uniquely has the form $P_{1}^{r}$ with 
$P_{1}$ of length divisible by $k$ and $r$ maximal. Then $r$ 
necessarily divides $n$ and $j$, and the orbit of $P$ under $R$ has 
size $n/r$.
In case $r\ge 1$, everything will merely be cut 
down by a factor of $r$. Declare two elements $(P_{1},i_{1})$ and 
$(P_{2},i_{2})$  to be \emph{rotation-equivalent} if $P_{1}$ and 
$P_{2}$ are in the same rotation class under $R$ (regardless of 
$i_{1}$ and $i_{2}$). 
As before, all elements of a rotation-equivalence class can be seen 
in a single diagram as illustrated.

\vspace*{-35mm}

\begin{pspicture}(-4,-2.5)(13,9)
\psset{unit=.5cm}

\multiput(0,9)(1,0){13}{.}
\multiput(0,8)(1,0){13}{.}
\multiput(0,7)(1,0){13}{.}
\multiput(0,6)(1,0){13}{.}
\multiput(0,5)(1,0){13}{.}
\multiput(0,4)(1,0){13}{.}
\multiput(0,3)(1,0){13}{.}
\multiput(0,2)(1,0){13}{.}
\multiput(0,1)(1,0){13}{.}
\multiput(0,0)(1,0){13}{.}

\psdots(0,0)
(1,1)(2,2)(3,3)(4,4)(5,5)(6,4)(7,5)(8,6)(9,7)(10,8)(11,9)(12,8)

\psdots(0,2)(2,4)(4,6)(6,2)(8,4)(10,6)

\psset{dotscale=2.3}
\blueclr{\psdots(0,0)(2,2)(4,4)(6,4)(8,6)(10,8)(12,8) }

\psline[linecolor=blue]
(0,0)(1,1)(2,2)(3,3)(4,4)(5,5)(6,4)(7,5)(8,6)(9,7)(10,8)(11,9)(12,8)

\psline(0,0)(6,2)
\psline(0,2)(6,4)
\psline(2,2)(8,4)
\psline(2,4)(8,6)
\psline(4,4)(10,6)
\psline(4,6)(10,8)

\psline[linestyle=dashed]{<-}(0,-1)(2,-1)
\psline[linestyle=dashed]{->}(4,-1)(6,-1)
\rput(3,-1){\textrm{{\small $P$}}}

\rput(6.5,2.1){\textrm{{\footnotesize $2$}}}
\rput(6.5,4.0){\textrm{{\footnotesize $1$}}}
\rput(8.5,4.2){\textrm{{\footnotesize $2$}}}
\rput(8.5,6.0){\textrm{{\footnotesize $0$}}}
\rput(10.5,6.1){\textrm{{\footnotesize $1$}}}
\rput(10.5,8.0){\textrm{{\footnotesize $0$}}}

\rput(6,-3.2){\textrm{{\small The rotation-equivalence class in
$\p_{3,2,2}^{*}$ of
$P=U^{5}D$}}}

\end{pspicture}

\vspace*{-35mm}

\begin{pspicture}(-1.5,-3.5)(25,8)
\psset{unit=.5cm}

\multiput(0,7)(1,0){25}{.}
\multiput(0,6)(1,0){25}{.}
\multiput(0,5)(1,0){25}{.}
\multiput(0,4)(1,0){25}{.}
\multiput(0,3)(1,0){25}{.}
\multiput(0,2)(1,0){25}{.}
\multiput(0,1)(1,0){25}{.}
\multiput(0,0)(1,0){25}{.}
\multiput(0,-1)(1,0){25}{.}
\multiput(0,-2)(1,0){25}{.}

\psdots(0,0)
(1,1)(2,2)(3,3)(4,2)(5,1)(6,0)(7,-1)(8,0)(9,1)(10,2)(11,1)(12,2)(13,3)(14,4)(15,5)
(16,4)(17,3)(18,2)(19,1)(20,2)(21,3)(22,4)(23,3)(24,4)

\psdots(0,2)(0,4)(3,5)(3,7)(6,2)(6,4)(9,3)(9,5)(12,-2)(12,0)(15,1)(15,3)(18,-2)(18,0)(21,-1)(21,1)

\psset{dotscale=2.3}
\blueclr{\psdots(0,0)(3,3)(6,0)(9,1)(12,2)(15,5)(18,2)(21,3)(24,4)}

\psline[linecolor=blue]
(0,0)(1,1)(2,2)(3,3)(4,2)(5,1)(6,0)(7,-1)(8,0)(9,1)(10,2)(11,1)(12,2)

\psline[linecolor=blue]
(12,2)(13,3)(14,4)(15,5)(16,4)(17,3)(18,2)(19,1)(20,2)(21,3)(22,4)(23,3)(24,4)

\psline(0,0)(12,-2)
\psline(0,2)(12,0)
\psline(0,4)(12,2)
\psline(3,3)(15,1)
\psline(3,5)(15,3)
\psline(3,7)(15,5)
\psline(6,0)(18,-2)
\psline(6,2)(18,0)
\psline(6,4)(18,2)
\psline(9,1)(21,-1)
\psline(9,3)(21,1)
\psline(9,5)(21,3)

\psline[linestyle=dashed]{<-}(0,-2.8)(5,-2.8)
\psline[linestyle=dashed]{->}(7,-2.8)(12,-2.8)
\rput(6,-2.8){\textrm{{\small $P$}}}

\rput(15.6,4.9){\textrm{{\footnotesize $0$}}}
\rput(15.5,2.9){\textrm{{\footnotesize $0$}}}
\rput(15.5,0.9){\textrm{{\footnotesize $1$}}}
\rput(12.5,1.9){\textrm{{\footnotesize $0$}}}
\rput(12.5,-0.1){\textrm{{\footnotesize $2$}}}
\rput(12.5,-2.1){\textrm{{\footnotesize $3$}}}
\rput(18.5,-2.1){\textrm{{\footnotesize $3$}}}
\rput(18.5,-0.1){\textrm{{\footnotesize $2$}}}
\rput(18.6,1.9){\textrm{{\footnotesize $1$}}}
\rput(21.5,-1.1){\textrm{{\footnotesize $3$}}}
\rput(21.5,0.9){\textrm{{\footnotesize $2$}}}
\rput(21.5,2.9){\textrm{{\footnotesize $1$}}}

\rput(12,-5){\textrm{{\small The rotation-equivalence class in
$\p_{4,3,3}^{*}$ of
$P=U^{3}D^{4}U^{3}DU$}}}

\end{pspicture}
Label the baselines (there are $jn/r$ of them; both illustrations have 
$r=1$) at their endpoints as follows (each of $0,1,\ldots,n-1$ will be the label on 
$j/r$ endpoints). First take the highest endpoint $p$ and consider the 
set of all endpoints lying weakly to the left of the vertical line 
through $p$. Since there are $j-1$ endpoints directly below $p$, 
this set has size at least $j$. Place label 0 on the $j/r$ highest 
points in this set, favoring points to the left if a choice must be 
made between points at the same height. Then take the highest 
unlabeled endpoint, consider the 
set of all unlabeled endpoints lying weakly to the left of its 
vertical line, 
and place the label 1 on the $j/r$ highest points in this 
set, again favoring ``left''. Continue in like manner until all endpoints 
are labeled.

Then, for each $i=0,1,\ldots,n-1$,the $j/r$ objects in the rotation-equivalence class with label 
$i$ all have $X=i$, and the uniform 
distribution of $X$ follows. By considering the objects in 
$\p_{n,k,j}^{*}$ with $X=n-1$, we obtain our main result.
\begin{main} Suppose $j\ge 1,\ k\ge 2,$ and $n\ge\frac{j}{k-1}$. Then 
$\frac{j}{n}\binom{kn}{n+j}$ is the number of lattice paths of $n+j$ 
upsteps $(1,1)$ and $kn-(n+j)$ downsteps $(1,-1)$ which $($i$\,)$ start at 
$(0,-2i)$ for some $i\ge 0$, and $($ii$\,)$ have all interior \kd points 
(strictly) above the line through the origin of slope 
$-\frac{(k-2)n-2}{kn}$.

\end{main} 

\vspace*{5mm}

\section{Concluding Remarks}

The main theorem can be generalized somewhat further (essentially the 
same proof): $\frac{d}{n}\binom{an}{cn+d}$ is the number of lattice 
paths of $cn+d$ upsteps and $an-(cn+d)$ downsteps which (i) start at 
$(0,-2i)$ for some $i\ge 0$, and (ii) have all interior \kd points 
(strictly) above the line through the origin of slope 
$-\frac{(a-2c)n-2}{an}$.

There is also a well known generalization of the Catalan numbers in a 
different direction: $\frac{j}{kn+j}\binom{kn+j}{n}$ is the number of 
lattice paths of $n$ steps east (1,0) and $(k-1)n+j-1$ steps north 
$(0,1)$ that start at the origin and lie weakly above the line 
$y=(k-1)x$. One way to prove this (slightly generalizing the approach 
in \cite{woan}) is as follows. Consider the set $\p_{n,k,j}$ of paths 
consisting of $n$ steps east and $(k-1)n+j$ steps north. Measuring 
``height'' of a point above $y=(k-1)x$ as the perpendicular distance 
to $y=(k-1)x$, define $j$ \emph{high points} for a path $P\in 
\p_{n,k,j}$: the first high point is the leftmost of the highest 
points on the path, the second high point is the leftmost of the next 
highest points of the path, and so on. Note that all $j$ high points 
necessarily lie strictly above $y=(k-1)x$. Mark any one of the these 
high points to obtain the set $\p_{n,k,j}^{*}$ of marked 
$\p_{n,k,j}$-paths. Clearly, $\v \p_{n,k,j}^{*} \v = 
j\binom{kn+j}{n}$. Label the $kn+j+1$ points on a marked path $P^{*} 
\in  \p_{n,k,j}^{*}$ in order $0,1,2,\ldots,kn+j$ starting at the 
origin. Set $X=$ label of the marked high point. Then $X$ is 
uniformly distributed over $1,2,\ldots,kn+j$. The paths with $X=kn+j$ 
yield the desired paths by deleting the last step (necessarily a north 
step) and rotating $180^{\circ}$.

All the above generalizations of the Catalan numbers are incorporated 
in the expression 
\[
\frac{ad-bc}{an+b}\binom{an+b}{cn+d}=(a-c)\binom{an+b-1}{cn+d-1}-c\binom{an+b-1}{cn+d}
\]
and it would be interesting to find a unified combinatorial 
interpretation for it.

\noindent 2000 {\it Mathematics Subject Classification}: 05A15.

\noindent \emph{Keywords: } Catalan, uniformly distributed, 
$k$-divisible, baseline, Cycle Lemma.


\begin{thebibliography}{99}
\bibitem{gkp} R.\,E.\,Graham, D.\,E.\,Knuth, Oren Patashnik, \emph{Concrete Mathematics} 
(2nd edition), Addison-Wesley, 1994.

\bibitem{zaks} N. Dershowitz and S. Zaks, The cycle lemma 
and some applications, \emph{European J. of Comb.} {\bf 11}, 1990, 35--40. 

\bibitem{sw} H. S. Snevily and D. B. West, The bricklayer 
problem and the strong cycle lemma,
\emph{Amer. Math. Monthly} {\bf105}, 1998,  131--143.


\bibitem{woan} Wen-jin Woan,
Uniform partitions of lattice paths and Chung-Feller generalizations. 
\emph{Amer. Math. Monthly} 108 (2001), no. 6, 556--559.

\end{thebibliography}
\end{document}